\documentclass{amsart}

\usepackage{amsmath,amsfonts}

\usepackage{amscd,amsthm}

\usepackage[dvips]{epsfig}

\newcommand{\R}{\mathop{\mathbb{R}}}

\newcommand{\Lip}[1]{\mathop{\mathit{Lip}\left( #1 \right)}}

\theoremstyle{plain} 
\newtheorem{theorem}{Theorem}[section]

\newtheorem{remark}[theorem]{Remark}

\numberwithin{equation}{section}

\begin{document}

\title{Uniqueness of limit cycles for a 
class of planar vector fields.}

\author{Timoteo Carletti}

\date{\today}

\address[Timoteo Carletti]{Scuola Normale Superiore, piazza dei Cavalieri 7,
 56126  Pisa, Italy}

\email[Timoteo Carletti]{t.carletti@sns.it}


\begin{abstract}
In this paper we give sufficient conditions to ensure uniqueness of limit cycles
for a class of planar vector fields. We also exhibit a class of examples with exactly one
limit cycle.
\end{abstract}

\maketitle

\section{Introduction}
\label{sec:intro}

In this paper we consider the problem of determine the number of limit cycles,
i.e. isolated closed trajectories, for planar vector fields. This is a classical
problem, included as part of the XVI Hilbert's problem. The literature is huge
and still growing, for a review see~\cite{Sansone} and the more recent~\cite{Zhang,Ye}.

An important subproblem is to study systems with a unique limit cycle, in fact
in this case the dynamics of the cycle can ''dominate'' the global dynamics of the
whole system.

Let us consider planar vector fields of the form:
\begin{equation}
\label{eq:vf}
\begin{cases}
\dot x &= \beta(x)\left[ \phi(y)-F(x,y)\right] \\
\dot y &= -\alpha(y)g(x) \, ,
\end{cases}
\end{equation}
under the regularity assumptions (to ensure the existence and
uniqueness of 
the Cauchy initial problem) there exists: $-\infty \leq a<0
<b \leq +\infty$, such that:
\begin{enumerate}
\item[A1)] $\beta\in \Lip{a,b}$ and $\alpha \in \Lip{\R}$;
\item[A2)] $\phi \in \Lip{\R}$, $g\in \Lip{a,b}$ and 
$F\in \Lip{(a,b)\times \R}$.
\end{enumerate}

Without loss of generality we can assume $\alpha$ and $\beta$ to be positive in
their respective domains of definition, in fact the existence of $x_0$ such that
$\beta(x_0)=0$ (or $y_0$ s.t. $\alpha(y_0)=0$), gives rise to 
{\em invariant lines}, which cannot intersect a limit cycle. Hence we can 
reparametrize time, by dividing the vector field by: $\alpha(y)\beta(x)$. 
The transformed system is:
\begin{equation}
\label{eq:vecfiel}
\begin{cases}
\dot x &= \tilde{\phi}(y)-\tilde{F}(x,y) \\
\dot y &= -\tilde{g}(x) \, ,
\end{cases}
\end{equation}
where $\tilde{\phi}(y)=\phi(y)/\alpha(y)$, $\tilde{F}(x,y)=F(x,y)/\alpha(y)$ and
$\tilde{g}(x)=g(x)/\beta(x)$. In the following we will drop out the 
$\tilde \,$--mark and consider the general system of previous type.

These systems can be though as ''non--Hamiltonian perturbations'' of Hamiltonian
 ones, with Hamilton function: $H(x,y)=\Phi(y)+G(x)$, where $\Phi(y)=
\int_0^y \phi(s)\, ds$ and $G(x)=\int_0^xg(s)\, ds$, being $F(x,y)$ the
''perturbation''. 

One can also consider~\eqref{eq:vecfiel} as ''generalized''
Li\'enard equations:
\begin{equation*}
\ddot x+f(x)\dot x+g(x)=0 \, ,
\end{equation*}
 which in the Li\'enard plane can be rewritten as:
\begin{equation}
\label{eq:lien}
\begin{cases}
\dot x &= y-F(x) \\
\dot y &= -g(x) \, ,
\end{cases}
\end{equation}
where $F^{\prime}(x)=f(x)$, hence our systems
 generalize~\eqref{eq:lien} by allowing a dependence of $F$ also on $y$.

Let $\lambda >0$ and let us consider the energy level $\mathcal{H}_{\lambda}=
\{ (x,y)\in \R^2: \Phi(y)+G(x)~=~\lambda \}$, the knowledge 
of the flow through $\mathcal{H}_{\lambda}$ can give informations about the
 existence of limit cycles. Because
\begin{equation*}
<\nabla \mathcal{H}_{\lambda},
X(x,y)>\Big \rvert_{\mathcal{H}_{\lambda}}=-F(x,y)g(x) \, ,
\end{equation*}
 where $X(x,y)= (\phi(y)-F(x,y),-g(x))$, no limit cycles can 
be completely contained in a 
region where $gF$ doesn't change sign. We will see in a while that the set
of zeros of $F$ will play a fundamental role in our construction.

For Li\'enard systems the set of zeros of $F$ is given by vertical lines
$x=x_k$ s.t. $F(x_k)=0$. In a recent 
paper~\cite{SabatiniVillari} authors, using ideas taken from Li\'enard 
systems~\cite{CarlettiVillari}, proved a uniqueness result for 
systems~\eqref{eq:vecfiel} assuming that $F(x,y)$ vanishes only at three vertical
lines $x=x_-<0$, $x=0$ and $x=x_+>0$. We generalize this condition by assuming that
zeros of $F(x,y)$ lie on (quite) general curves. More precisely let us assume 
there exist 
 $\psi_j:\R \rightarrow \R$, $j\in \{1,2\}$, such 
that~\footnote{We remark that our main result
still holds, even if one assume there exist 
$\alpha_j<0<\beta_j$, $j\in \{1,2\}$,
 and the functions $\psi_j$ to be defined in $[\alpha_j,\beta_j]$ and verify
hypothesis B) on their new domain of definition.}:
\begin{enumerate}
\item[B0)] $F(0,y)=0$ for all real $y$;
\item[B1)] $y\mapsto\psi_1(y)$, is {\em positive} for all $y\in \R $,
increasing for negative $y$, decreasing for positive $y$, $\psi_1(0) < b$;
\item[B2)] $y\mapsto\psi_2(y)$, is {\em negative} for all $y\in \R$,
decreasing for negative $y$, increasing for positive $y$, $\psi_2(0) > a$;
\item[B3)] for all $y\in \R$, $j\in \{ 1,2 \}$, we have:
%
\begin{equation*}
F(\psi_j(y),y)\equiv 0 \, ,
\end{equation*}
these curves will be called ''non--trivial zeros'' of $F(x,y)$ (in opposition with
the trivial zeros given by $x=0$).
\end{enumerate}

Let us divide the strip $(a,b)\times \R $
 into four distinct domains:
\begin{itemize}
\item $D_1^{>}:=\{ (x,y) \in (a,b)\times \R : x > \psi_1(y) \}$;
\item $D_1^{<}:=\{ (x,y) \in (a,b)\times \R : 0<x < \psi_1(y) \}$;
\item $D_2^{>}:=\{ (x,y) \in (a,b)\times \R : \psi_2(y) < x < 0 \}$;
\item $D_2^{<}:=\{ (x,y) \in (a,b)\times \R : x < \psi_2(y) \}$.
\end{itemize}

The following assumptions generalize ''standard sign ones'':
\begin{enumerate}
\item[C1)] $y\phi(y)>0$ for all $y\neq 0$ and $xg(x) >0$ for all 
$x\in (a,b)\setminus \{ 0 \}$; 
\item[C2)] $g(x)F(x,y)<0$ for all $(x,y)\in D_1^< \cup D_2^>$.
\end{enumerate}

We remark that hypothesis C2) can be weakened into:
\begin{enumerate}
\item[C2')] $g(x)F(x,y)\leq 0$ for all $(x,y)\in D_1^< \cup D_2^>$
 except at some $(x_0,y_0)$ where strictly inequality holds.
\end{enumerate}

With these hypotheses we ensures that $(0,0)$ is the only singular
point of system~\eqref{eq:vecfiel} and trajectories wind clockwise around it.

We are now able to state our main result
\begin{theorem}
\label{thm:main}
Let us consider system~\eqref{eq:vecfiel} and let us assume Hypotheses A), B) 
and C) to hold. Then there is at most one limit cycles which intersect both curves
$x=\psi_1(y)$ and $x=\psi_2(y)$ contained in $(a,b)\times \R$, provided:
\begin{enumerate}
\item[D1)] the function $y\mapsto F(x,y)/\phi(y)$ is strictly increasing
for $(x,y) \in D_1^{<}$ and $y\neq 0$;
\item[D2)] the function $y\mapsto F(x,y)/\phi(y)$ is strictly decreasing
for $(x,y) \in D_2^{>}$ and $y\neq 0$;
\item[E)] the function $x\mapsto F(x,y)$ is positive in $D_1^{>}$, negative
in $D_2^{<}$ and increasing in $D_1^{>} \cup D_2^{<}$;
\item[F)] let $A_j(y)=\left[\phi(y)\partial_x F(x,y)-g(x)\partial_y 
F(x,y)\right]\Big \rvert_{x=\psi_j(y)}$, then $A_j(y)y>0$ for $y\neq 0$, $j\in \{1,2\}$.
\end{enumerate}

\end{theorem}
Hypotheses D) and E) naturally generalize hypotheses used in the Li\'enard 
case~\cite{CarlettiVillari} or in the more general situation studied 
in~\cite{SabatiniVillari}. Also hypothesis F) is very natural: each closed 
trajectory intersects the non--trivial zeros of $F(x,y)$ at most once in any
 quadrant. We remark that this condition is trivially verified if the functions
$\psi_j$ are indentically constant, namely in the case considered 
in~\cite{SabatiniVillari}.

The proof of this result will be given in the next section. In the last
section (\S~\ref{sec:example}) we will provide a family of systems with
 exactly one limit cycle.

\section{Proof of Theorem~\ref{thm:main}}
\label{sec:proof}

The aim of this section is to prove our main result Theorem~\ref{thm:main}. 
The proof is based on the following remark,
\begin{remark}
\label{rem:integral}
Along any closed curve $\gamma:[0,T]\rightarrow (a,b)\times \R$ one has:
\begin{equation*}
\int_0^T \frac{d}{dt}H\Big \rvert_{\text{flow}}\circ \gamma(s) \, ds = H\circ \gamma (T)- H\circ \gamma (0)
=0 \, ,
\end{equation*}
moreover if $\gamma$ is an integral curve of system~\eqref{eq:vecfiel} we can 
evaluate the integrand function to obtain:
\begin{equation}
\label{eq:intzero}
I_{\gamma}:=\int_0^T g(x_{\gamma}(s)) F(x_{\gamma}(s),y_{\gamma}(s)) \, ds = 0 \, ,
\end{equation}
where $\gamma(s)=(x_{\gamma}(s),y_{\gamma}(s))$.
\end{remark}
The uniqueness result will be proved by showing that the existence of two
limit cycles, $\gamma_1$ contained~\footnote{By this we mean $\gamma_1$
is properly contained in the compact set whose boundary is $\gamma_2$.}
 in $\gamma_2$, both intersecting
$x=\psi_1(y)$ and $x=\psi_2(y)$, will imply: 
$I_{\gamma_1}< I_{\gamma_2}$, which contradicts~\eqref{eq:intzero}.

From now on we will assume the existence of two limit cycles,
 $\gamma_1$ contained
in $\gamma_2$, which intersect both non--trivial zeros of $F$.

Let us now consider the set of zeros of the equation: $\phi(y)-F(x,y)=0$ inside
$D_1^< \cup D_2^>$; we claim that this set is the graph of some function $x\mapsto
\zeta(x)$. Moreover this function vanishes for $x\in \{\psi_2(0),0,\psi_1(0) \}$,
it is positive for $x\in (\psi_2(0),0)$ and negative for $x\in (0,\psi_1(0))$.

The first statement follows from hypotheses D). In fact, let assume the claim to be 
false, then there exist $(\bar x,y_k)\in D_2^>$, such that:
\begin{equation*}
y_1 < y_2 \quad \text{and} \quad \phi(y_k)-F(\bar x,y_k)=0 \quad k\in\{1,2\} \, .
\end{equation*}
Then using the decreasing property (hypothesis D)) we obtain a contradiction:
\begin{equation*}
0=\left[F(\bar x,y_1)-\phi(y_1)\right]/\phi(y_1)>
\left[F(\bar x,y_2)-\phi(y_2)\right]/\phi(y_2)=0\, .
\end{equation*}
We can do the same construction for $D_1^<$, and the statement is proved.
The sign properties of $\zeta$ can be proved as follows. Let
$\psi_2(0)<x<0$, 
then by C2) $0<F(x,\zeta(x))=\phi(\zeta(x))$,
 using now C1) we conclude that $\zeta(x)>0$. The case $0<x<\psi_1(0)$
 can be handle similarly and we omit. By continuity we get the result
 about the zeros of $\zeta(x)$.

Hypothesis F) guarantees that a closed trajectory can intersect the non--trivial
zeros of $F(x,y)$ only once in each quadrant, in fact $A_j(y)$ gives a measure of 
the angle between the vector field and the normal to 
$\mathcal{F}_0=\{ (x,y): x=\psi_1(y) \} \cup \{ (x,y): x=\psi_2(y) \}$
 at $(\psi_j(y),y)$:
\begin{equation*}
<\nabla \mathcal{F}_{0},
X(x,y)\Big \rvert_{\mathcal{F}_{0}}>=\left[\phi(y)\partial_x F(x,y)
-g(x)\partial_y F(x,y)\right]\Big\rvert_{\mathcal{F}_0} =A_j(y)\quad j\in
\{1,2 \} \, .
\end{equation*}
For instance, because the angle between the vector field and 
$\{ (x,y): x=\psi_1(y) \}\cap \{ y >0 \}$ is in absolute value smaller than $\pi/2$,
 a trajectory starting at $(0,\bar{y})$, for some $\bar{y}>0$, which will
intersect $\{ (x,y): x=\psi_1(y) \}\cap \{ y >0 \}$, could not meet anew 
$\{ (x,y): x=\psi_1(y) \}\cap \{ y >0 \}$.
 
From hypotheses D) and the sign of $F$ on $D_2^< \cup D_1^>$, it follows easily
 that for all $(x,y)\in D_2^< \cup D_1^>$ one has: 
$\left(y-\zeta(x)\right)\left(\phi(y)-F(x,y)\right)>0$. Hence a cycle intersects
$\left(D_2^< \cup D_1^>\right)\cap \{ y>0 \}$ in a region where $\phi(y)-F(x,y)>0$,
whereas the intersection with $\left(D_2^< \cup D_1^>\right)\cap \{ y<0 \}$ holds
where $\phi(y)-F(x,y)<0$. This remark allows us to divide the path of integration
needed to evaluate $I_{\gamma_j}$, $j\in \{1,2\}$, in two parts: an ''horizontal''
one where $\dot x >0$ and a vertical one,
 where $\dot x$ vanishes.

\begin{center}
  \begin{figure}[ht]
   \makebox{
   \epsfig{file=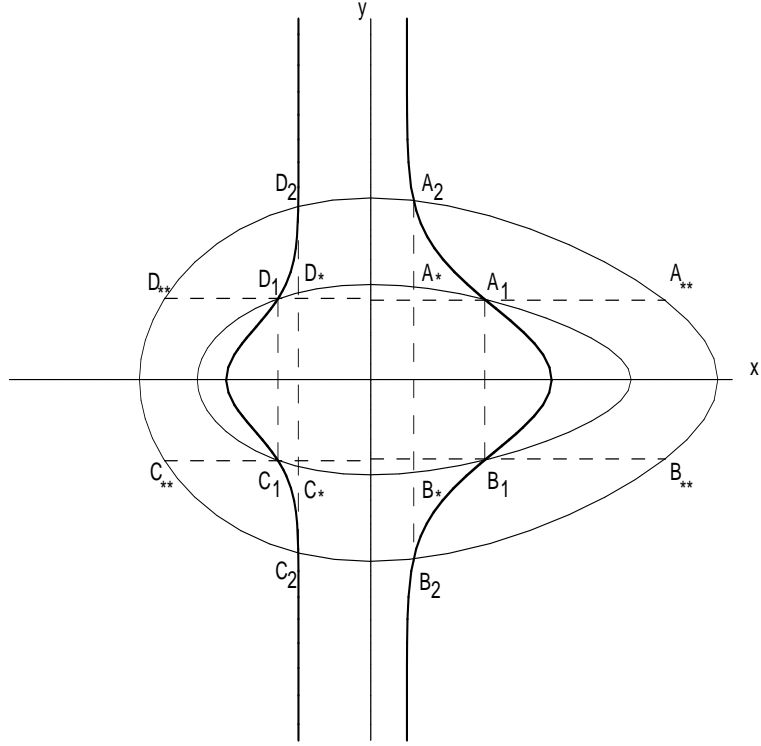,height=12cm,width=12cm,angle=-90}}
  \caption{The non--trivial zeros of $F$ (thick), the limit cycles $\gamma_1$
and $\gamma_2$ (thin) intersecting both non--trivial zeros of $F$ 
and their subdivision into arcs.}
  \label{fig:figura1}
  \end{figure}
\end{center}

Let us define (see Figure~\ref{fig:figura1}), for $j\in \{1,2\}$, $A_j$ 
(respectively $B_j$) the intersection 
point of $\gamma_j$ with $x=\psi_1(y)$ for $y>0$ (respectively $y<0$), and 
$C_j$ (respectively $D_j$) the intersection 
point of $\gamma_j$ with $x=\psi_2(y)$ for $y<0$ (respectively $y>0$).
Let also introduce, $A_*$ being the intersection point of $\gamma_1$ and the
line $x=x_{A_2}$ contained in the first quadrant, and $A_{**}$ being the
 intersection point of $\gamma_2$ and the
line $y=y_{A_1}$ contained in the first quadrant. Similarly we introduce points:
$B_*$, $B_{**}$, $C_*$, $C_{**}$  and $D_*$, $D_{**}$ 
(see Figure~\ref{fig:figura1}).

According to this subdivision of the arcs of limit cycles, we evaluate 
$I_{\gamma_j}$ as follows:
\begin{eqnarray}
\label{eq:sudiv}
I_{\gamma_1}=\int_{ {D_*A_*}}+\int_{{A_*A_1}}+
\int_{ {A_1B_1}}+
\int_{ {B_1B_*}}+
\int_{ {A_*C_*}}+
\int_{ {C_*C_1}}+
\int_{ {C_1D_1}}+
\int_{ {D_1D_*}}  \\
I_{\gamma_2}=\int_{ {D_2A_2}}+\int_{ {A_2A_{**}}}
+\int_{ {A_{**}B_{**}}}+\int_{ {B_{**}B_2}}
+\int_{ {B_2C_2}}+\int_{ {C_{2}C_{**}}}
+\int_{ {C_{**}D_{**}}}+\int_{ {D_{**}D_2}} \notag\, .
\end{eqnarray}

Let now show that $I_{\gamma_1}<I_{\gamma_2}$, which prove the contradiction
and conclude the proof.

\subsection{Integration along ''horizontal arcs''.}
\label{ssec:intoriz}

Because along horizontal arcs we have $\dot x \neq 0$, we can change integration
variable from $t$ to $x$, hence for example:
\begin{equation*}
\int_{ {D_jA_j}}g(x)F(x,y)\, dt= \int_{x_{D_j}}^{x_{A_j}} \frac{g(x)
F(x,y_j(x))}{\phi(y_j(x))-F(x,y_j(x))} \, dx \, ,
\end{equation*}
where $y_j(x)$, $j\in \{1,2 \}$, is the parametrization of $\gamma_j$ as 
graph over $x$ for $x\in(x_{D_j},x_{A_j})$.

Because $y_2(x)>y_1(x)$ for all $x\in (x_{D_*},x_{A_*})$, using hypotheses D) 
and the sign assumptions C) we get: 
\begin{equation*}
\frac{g(x)F(x,y_1(x))}{\phi(y_1(x))-F(x,y_1(x))} <
\frac{g(x)F(x,y_2(x))}{\phi(y_2(x))-F(x,y_2(x))} \, ,
\end{equation*}
hence:
\begin{eqnarray}
\label{eq:d1a1}
\int_{x_{D_{1}}}^{x_{A_{1}}} \frac{g(x)
F(x,y_1(x))}{\phi(y_1(x))-F(x,y_1(x))} \, dx & < 
\int_{x_{D_{1}}}^{x_{D_{*}}} \frac{g(x)
F(x,y_1(x))}{\phi(y_1(x))-F(x,y_1(x))} \, dx+\int_{x_{A_*}}^{x_{A_1}} \frac{g(x)
F(x,y_1(x))}{\phi(y_1(x))-F(x,y_1(x))} \, dx \notag \\ 
&+\int_{x_{D_2}}^{x_{A_2}} \frac{g(x) F(x,y_2(x))}{\phi(y_2(x))-F(x,y_2(x))} \, dx 
\leq \int_{x_{D_2}}^{x_{A_2}} \frac{g(x)
F(x,y_2(x))}{\phi(y_2(x))-F(x,y_2(x))} \, dx \, ,
\end{eqnarray}
the last step follows because the integrand function is negative by hypothesis C2)
and from the previous discussion on the sign of $\phi(y)-F(x,y)$.

In a very similar way we can prove that:
\begin{equation}
\label{eq:b1c1}
\int_{x_{B_{1}}}^{x_{C_{1}}} \frac{g(x)
F(x,y_1(x))}{\phi(y_1(x))-F(x,y_1(x))} \, dx <
\int_{x_{B_2}}^{x_{C_2}} \frac{g(x) F(x,y_2(x))}{\phi(y_2(x))-F(x,y_2(x))} \, dx 
\, .
\end{equation}

\subsection{Integration along ''vertical arcs''.}
\label{ssec:intvert}

Along vertical arcs $\dot y$ never vanishes, hence we can perform the
integration w.r.t. the $y$ variable and getting for example:
\begin{equation*}
\int_{ {A_jB_j}}g(x)F(x,y)\, dt= \int_{y_{B_j}}^{y_{A_j}} F(x_j(y),y)\, dy
 \, ,
\end{equation*}
where $x_j(y)$, $j\in \{1,2 \}$, is the parametrization of $\gamma_j$ as 
graph over $y$ for $y\in(y_{B_j},y_{A_j})$.

Because $x_2(y)>x_1(y)$ for all $y\in (y_{A_{**}},y_{B_{**}})$, from hypotheses E) 
and the definition of $A_{**}$ and $B_{**}$, we get:
\begin{equation*}
\int_{y_{B_{1}}}^{y_{A_{1}}} F(x_1(y),y)\, dy < 
\int_{y_{B_{**}}}^{y_{A_{**}}} F(x_2(y),y)\, dy  \, .
\end{equation*}
Again from the sign assumption on $F$ in $D_1^{>}$, we get:
\begin{equation*}
\int_{y_{A_{**}}}^{y_{A_{2}}} F(x_2(y),y)\, dy  >0 \quad\text{and}\quad
\int_{y_{B_{2}}}^{y_{B_{**}}} F(x_2(y),y)\, dy >0 \, ,
\end{equation*}
hence we obtain:
\begin{equation}
\label{eq:a2b2}
\int_{y_{B_1}}^{y_{A_1}} F(x_1(y),y)\, dy 
<\int_{y_{B_2}}^{y_{A_2}} F(x_2(y),y)\, dy \, .
\end{equation}
Analogously we can prove that:
\begin{equation}
\label{eq:c2d2}
\int_{y_{C_1}}^{y_{D_1}} F(x_1(y),y)\, dy 
<\int_{y_{C_2}}^{y_{D_2}} F(x_2(y),y)\, dy \, .
\end{equation}

\subsection{Conclusion of the proof.}
\label{ssec:concpro}

We are now able to complete our proof. In fact from~\eqref{eq:d1a1} 
and~\eqref{eq:b1c1} of \S~\ref{ssec:intoriz}, from~\eqref{eq:a2b2} 
and~\eqref{eq:c2d2} of \S~\ref{ssec:intvert} and the subdivision~\eqref{eq:sudiv}
we get:
\begin{equation*}
I_{\gamma_1}<I_{\gamma_2} \, ,
\end{equation*}
which contradicts~\eqref{eq:intzero}, and so the Theorem is proved.

\section{A system with exactly one limit cycle}
\label{sec:example}

In this last part we present an important class of examples exhibiting exactly
one limit cycle. Let us
assume that $F$ has the following ''special form'': 
\begin{equation*}
F(x,y)=x\left[x-\psi_1(y)\right]\left[x-\psi_2(y)\right] \, ,
\end{equation*}
where $\psi_j$ verify 
hypotheses B). Then hypothesis F) is equivalent to:
\begin{enumerate}
\item[F')] the function $y\mapsto \Phi(y)+G(\psi_j(y))$ is strictly increasing
for positive $y$ and strictly decreasing for negative ones, $j \in \{ 1,2 \}$.
\end{enumerate}
In fact we have:
\begin{eqnarray*}
A_1(y)&=\psi_1(y)\left[ \psi_1(y) - \psi_2(y) \right] \left[ \phi(y)-g(\psi_1(y)) 
\psi_1^{\prime}(y)\right]\\
&=\psi_1(y)\left[\psi_1(y)-\psi_2(y)\right]\frac{d}{dy}\left[
\Phi(y)-G(\psi_1(y))\right] \, ,
\end{eqnarray*}
and the claim follows from the sign properties of $\psi_j$ and the definitions of
$\Phi$ and $G$. Similarly for $A_2$.

A concrete example is given by choosing:
\begin{equation*}
\phi(y)=y \, , \, g(x)=x \, , \, \psi_1(y)=c_1 e^{-d_1 y^2}+e_1 \quad \text{and}
 \quad \psi_2(y)=-c_2 e^{-d_2 y^2}-e_2 \, ,
\end{equation*}
 with $c_j,d_j,$ and $e_j$ positive real
numbers such that:
\begin{enumerate}
\item $c_1+e_1=c_2+e_2$,
\item let $r=c_1+e_1$, then $c_jd_j\max\{ r, r^2\}<1/2$ for $j\in \{ 1,2 \}$.
\end{enumerate}

\begin{center}
  \begin{figure}[ht]
   \makebox{
   \epsfig{file=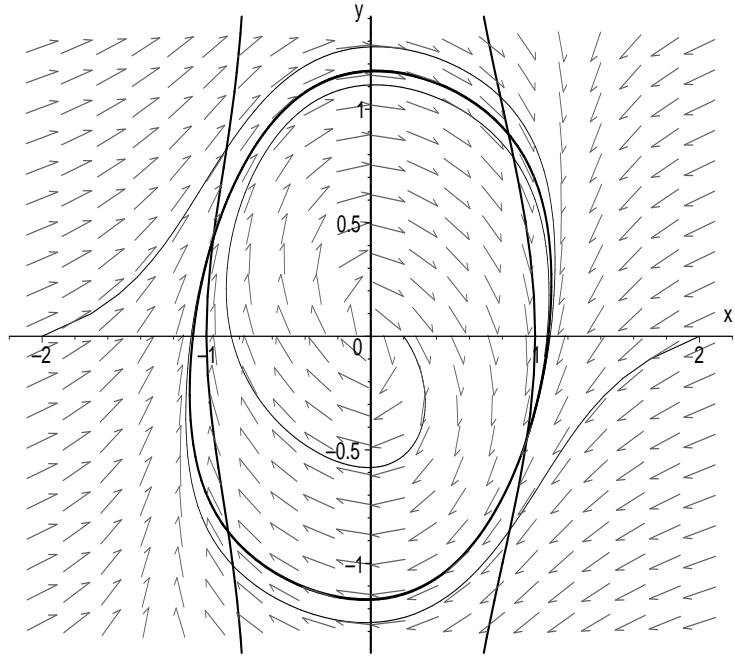,height=12cm,width=12cm,angle=-90}}
  \caption{An example with $c_1=d_1=e_1=1/2$, $c_2=1/4$, $d_2=1$, and $e_2=3/4$.
 We numerically compute the attracting limit cycle (thick) and three attracted
 trajectories (thin), we also plot the non--trivial zeros of $F$
 (thick) and the direction of
the vector field (arrows).}
  \label{fig:figura2}
  \end{figure}
\end{center}

The example of Figure~\ref{fig:figura2} verifies all the hypotheses of Theorem~\ref{thm:main}. 
Hypotheses A), B), C) and E) are readily verified, whereas hypotheses 
D) and F') (hence F)) follow~\footnote{Using $c_jd_jr<1/2$ hypothesis
F') holds. Whereas (1) and $c_jd_jr^2<1/2$ ensure that the circle $x^2+y^2=r^2$
is contained inside $D_1^< \cup D_2^>$ and it is tangent to $x=\psi_j(y)$ at
$y=0$, for $j\in \{ 1,2\}$, thus any limit cycle must intersect both
curves $x=\psi_j(y)$, $j\in \{ 1,2\}$.}
 from (1) and (2) plus some easy calculations that
we omit.

\end{document}